\documentclass[letterpaper, 10 pt, conference]{ieeeconf}  

\IEEEoverridecommandlockouts                              
\usepackage{amsmath} 
\usepackage{amssymb}  
\usepackage{xcolor}

\usepackage{algorithm} 
\usepackage{algpseudocode} 

\usepackage{booktabs}
\usepackage{array}
\usepackage{soul}
\usepackage{graphicx}

\newcommand{\allocation}{z}
\newcommand{\relaxation}{\text{rel}}
\newcommand{\utility}{u}

\title{\LARGE \bf
A Linear and Scalable Cutting-Plane Algorithm for Electricity Pricing 
}

 \author{ \parbox{2 in}{\centering Mat\'ias Romero
        \\
        GSB, Columbia University\\
        New York, NY 10027, USA\\
        {\tt\small mer2262@columbia.edu}}
        \hspace*{ 0.2 in}
        \parbox{2 in}{ \centering Felipe Ver\'astegui
        \\
       IEOR, Columbia University\\
        New York, NY 10027, USA\\
        {\tt\small fav2114@columbia.edu}}
        \hspace*{ 0.2 in}
        \parbox{2 in}{ \centering Mat\'ias Villagra
        \thanks{* Corresponding author: Matías Villagra, mjv2153@columbia.edu}
        \\
       IEOR, Columbia University\\
        New York, NY 10027, USA\\
        {\tt\small mjv2153@columbia.edu}}
}

\begin{document}

\maketitle
\thispagestyle{empty}
\pagestyle{empty}

\begin{abstract}

We propose a linear cutting-plane pricing algorithm tailored for large-scale electricity markets, addressing nonconvexities arising from the Alternating Current Optimal Power Flow equations. We benchmark our algorithm against a Direct Current (DC) approximation and the Jabr Second-Order Cone (SOC) relaxation under both the Integer Programming and Convex Hull pricing rules. We provide numerical results for a small (617-bus) and three large ($\geq 15,000$-bus) networks. Our algorithm yields price signals very close to the Jabr SOC, with computation times comparable to DC once we allow for warm-starts, including scenarios with line contingencies. 

\end{abstract}


\section{Introduction}

Price formation in electricity markets is 
inherently challenging. Feasible allocations for the Alternating Current Optimal Power Flow (AC-OPF) with unit commitments, belong to a nonlinear, nonconvex, and possibly disconnected region~\cite{hiskes_2001_exploring, bukhsh_2013_localsolutions}. One option is to use the Direct Current (DC) linear approximation which simplifies the pricing problem, but may yield solutions that are far from the AC-feasible region, requiring post-market clearing adjustments by Independent System Operators (ISOs)
to recover AC-feasibility. An alternative is to use tight convex nonlinear relaxations such as the Jabr Second-Order Cone (SOC) or the Quadratic Convex (QC) relaxation~\cite{jabr_2006_conic, coffrin_2016_qc}. These formulations yield less biased price signals on small grids~\cite{bichler_2023_prices}, but face scalability issues in medium to large grids~\cite{coffrin_2016_qc, bichler_2023_prices, villagra_2024_accurate}. An open research question is whether convex relaxations can be used for price formation in large grids, including reactive power prices, improving the signaling of congestion and scarcity in AC-feasible outcomes in market incentives~\cite{bichler_2023_prices}. 

A second source of nonconvexity originates from the sellers, where binary commitments over generators induce nonconvex costs. Much of the existing literature on electricity pricing focuses on rules to deal with these discrete costs, such as commitment and start-up costs, which pose additional challenges on top of the AC-OPF \cite{mays_2024_sequential}. Two well-known rules that address this issue are Convex Hull (CH) and Integer Programming (IP) pricing \cite{oneil_2005_ipprices, gribik_2007_marketclearing}. Both require that the underlying power flow model remains {\it convex}. The DC-OPF approximation is a special case where the underlying model is {\it linear}, making it a preferred companion choice for IP and CH pricing. The SOC and QC relaxations trade linearity for a more accurate representation of the AC power flow equations, but this comes at the cost of increased computational effort when combined with the CH or IP methodologies~\cite{bichler_2023_prices}.

In this paper, we bridge a gap in the electricity pricing literature by developing a scalable pricing algorithm which produces a very tight and numerically stable {\it linear} relaxation for AC-OPF, achieving computational performance comparable to the DC-OPF formulation, while producing active and reactive power prices which are very close to the Jabr SOC. To this end, we leverage the Cutting-Plane (CP) algorithm developed in~\cite{villagra_2024_accurate} which tightly outer-approximates the Jabr SOC. We run computational experiments under IP and CH pricing rules. In a small grid ($< 1,000$ nodes), we find that our algorithm with IP pricing yields allocations with very low redispatch costs compared to DC, and very accurate prices compared to the Jabr SOC, within a fraction of its computation time. In large grids ($\geq 15,000$ nodes), we find that the CP pricing algorithm with CH pricing scales similarly to DC when the cuts can be warm-started, and again the accuracy remains very high compared to the Jabr SOC.


The paper is organized as follows: Section~\ref{sec:modeling-framework} presents the modeling framework; Section~\ref{sec:algorithmic-price-formation} describes our price formation algorithm; Section~\ref{sec:numerical-results} presents computational experiments; and Section~\ref{sec:conclusions-future-work} concludes. 


\vspace{-2mm}
\section{Modeling Framework} \label{sec:modeling-framework}

Let $\mathcal{N} := (\mathcal{B},\mathcal{E})$ denote a power system where $\mathcal{B}$ denotes the set of buses and $\mathcal{E}$ denotes the set of branches. Price formation starts with a welfare problem, where buyers or loads, denoted by $l \in \mathcal{L}$ (sellers or generators, denoted by $g \in \mathcal{G}$) are assumed to truthfully submit piece-wise linear bidding functions $\pi_{\tau}: \mathbb{R}^{+} \mapsto \mathbb{R}^{+}$, $\tau \in \{l,g\}$ mapping consumption (production, binary commitments) to profits (costs). For each bus $k \in \mathcal{B}$, $\mathcal{G}_{k} \subseteq \mathcal{G}$ and $\mathcal{L}_{k} \subseteq \mathcal{L}$ denotes the generators and loads at bus $k$. The problem can be written as 

\begin{subequations}\label{model:acopf-welfare-problem}
{\small 
\begin{align}
&\max \hspace{1em} \sum_{l \in \mathcal{L}} \pi_{l}(p_l) + \sum_{g \in \mathcal{G}} \pi_{g} (p_g, x^{on}_g, x^{su}_g, x^{sd}_g)  \label{eq:acopf-welfare-obective} \\
     & \text{subject to:}  \hspace{5em} \nonumber\\
    & \sum_{km \in \delta(k)}  P_{km} = \sum_{g \in \mathcal{G}_{k}} p_{g} - \sum_{\ell \in \mathcal{L}_{k}} p_{\ell}  \quad \forall k \in \mathcal{B} \label{eq:acopf-welfare-active-p-balance} \\
   & \sum_{km \in \delta(k)} Q_{km}  = \sum_{g \in \mathcal{G}_{k}} q_{g} -\sum_{\ell \in \mathcal{L}_{k}} q_{\ell}   \quad \forall k \in \mathcal{B} \label{eq:acopf-welfare-reactive-p-balance} \\
    & (x^{on}, x^{su},  x^{sd}) \in \{0,1\}^{3 \times \mathcal{G}} \label{eq:acopf-welfare-binary-commit} \\
    & (p_g,  q_g)_{\forall g \in  \mathcal{G}} \in \mathcal{P}_{\mathcal{G}}(x^{on}, x^{su}, x^{sd}) \label{eq:acopf-welfare-gens-feasibility} \\
    & (p_l, q_l)_{\forall l \in \mathcal{L}}  \in \mathcal{P}_{\mathcal{L}} \label{eq:acopf-welfare-loads-feasibility} \\
    & (P_{km}, Q_{km})_{\forall km \in \epsilon} \in \mathcal{Q} \label{eq:acopf-welfare-power-flow-region}
\end{align}
}
\end{subequations}
\noindent where a social planner maximizes surplus subject to: active and reactive power balance \eqref{eq:acopf-welfare-active-p-balance}-\eqref{eq:acopf-welfare-reactive-p-balance}; binary commitments for generators \eqref{eq:acopf-welfare-binary-commit}; sellers~\eqref{eq:acopf-welfare-gens-feasibility} and buyers~\eqref{eq:acopf-welfare-loads-feasibility} production and consumption constraints; and a representation of the power flows in steady state and operational constraints~\eqref{eq:acopf-welfare-power-flow-region}. The set $\mathcal{P}_{\mathcal{L}}$ is polyhedral, and given a binary commitment $(x_{0}^{on}, x_{0}^{su}, x_{0}^{sd})$, the set $\mathcal{P}_{\mathcal{G}}(x_{0}^{on}, x_{0}^{su}, x_{0}^{sd})$ is polyhedral as well. The set $\mathcal{Q}$ is bounded and generally nonconvex.

\subsection{Full AC-OPF Formulation}
 A complete description of the physical constraints for active and reactive power flows in \eqref{model:acopf-welfare-problem} is given by the AC-OPF.  Using the {\it polar formulation}, the set $\mathcal{Q}$ will be denoted by $\mathcal{Q}^{\text{AC}}$ and can be described by:

\vspace{-\baselineskip}

{\small
\begin{subequations} \label{eq:Q_AC}
    \begin{align}
        & \left\{ \begin{aligned}
        V^{\min}_{k} \leq | V_{k} | \leq V^{\max}_{k} \nonumber
        \end{aligned} \right\}_{\forall k \in \mathcal{B}} \\
        & \left\{ \begin{aligned} 
        | \theta_{km} | &\leq \bar{\Delta}_{km}, \quad \theta_{km} = \theta_{k} - \theta_{m} \\
        P_{km} &= G_{kk} |V_k|^{2} + G_{km} |V_k| |V_m| \cos(\theta_{km}) \nonumber  \\
        &\hspace{1em} + B_{km}|V_k| |V_m| \sin(\theta_{km}) \\
        P_{mk} &= G_{mm} |V_m|^{2} + |V_k| |V_m| G_{mk}  \cos(\theta_{km}) \nonumber \\
        &\hspace{1em} - B_{mk}|V_k| |V_m| \sin(\theta_{km}) \\
        Q_{km} &= -B_{kk} |V_k|^{2} + B_{km} |V_k| |V_m| \cos(\theta_{km}) \nonumber \\
        &\hspace{1em} - G_{km}|V_k| |V_m| \sin(\theta_{km}) \\
        Q_{mk} &= -B_{mm} |V_m|^{2} + |V_k| |V_m| B_{mk} \cos(\theta_{km}) \nonumber \\
        &\hspace{1em} + G_{mk}|V_k| |V_m| \sin(\theta_{km}) \\
        P^{2}_{km} &+ Q^{2}_{km} \leq \overline{i}^{2}_{km} |V_{k}|^{2} \\
        P^{2}_{mk} &+ Q^{2}_{mk} \leq \overline{i}^{2}_{km} |V_{m}|^{2}  \\
        \end{aligned} \right\}_{\forall km \in \mathcal{E}} 
    \end{align}
\end{subequations}
}
where $\theta_{km}, \theta_k$ and $V_k$ are voltage angle and magnitude variables. The physical parameters of each branch $km \in \mathcal{E}$ are described by
\begin{equation*}
    Y_{km} := \begin{pmatrix}
    G_{kk} + j B_{kk} & G_{km} + j B_{km} \\
    G_{mk} + j B_{mk} & G_{mm} + j B_{mm},
    \end{pmatrix}
\end{equation*}
which is the complex admittance matrix of branch $km$, and $\overline{\Delta}_{km}$ denotes maximum angle difference and $\overline{i}^{2}_{km}$ a bound on apparent current-squared.

\subsection{Second-Order Cone Relaxations}

The use of $\mathcal{Q}^{\text{AC}}$ in \eqref{model:acopf-welfare-problem} poses a challenge for competitive equilibrium pricing since $\mathcal{Q}^{\text{AC}}$ is nonconvex~\cite{bichler_2023_prices, bickchandani_2002_assignment}. 
The power flow equations in $\mathcal{Q}^{\text{AC}}$ can be linearized using additional variables $v_{k}^{(2)}, c_{km}, s_{km}$, and by adding appropriate rotated-cone inequalities, known as \emph{Jabr} inequalities~\cite{jabr_2006_conic}, we obtain the convex set $\mathcal{Q}^{\text{Jabr}}$ defined by:

\vspace{-\baselineskip}

{\small
\begin{subequations}
    \begin{align}
        & \left\{ \begin{aligned}  \nonumber
        v^{(2), \min}_{k} \leq v^{(2)}_{k} \leq v^{(2), \max}_{k}
        \end{aligned} \right\}_{\forall k \in \mathcal{B}} \\
        & \left\{ \begin{aligned} \nonumber
        c^2_{km} &+ s^2_{km} \leq v^{(2)}_{k} v^{(2)}_{m} \\
        P_{km} &= G_{kk} v^{(2)}_{k} + G_{km} c_{km} + B_{km} s_{km} \\
        P_{mk} &= G_{mm} v^{(2)}_{m} + G_{mk} c_{km} - B_{mk} s_{km} \\
        Q_{km} &= -B_{kk} v^{(2)}_{k} + B_{km} c_{km} - G_{km} s_{km} \\
        Q_{mk} &= -B_{mm} v^{(2)}_{m} + B_{mk} c_{km} + G_{mk} s_{km} \\
        P^{2}_{km} &+ Q^{2}_{km} \leq \overline{i}^{2}_{km} v_{k}^{(2)} \\
        P^{2}_{mk} &+ Q^{2}_{mk} \leq \overline{i}^{2}_{km} v_{m}^{(2)}  \\
        \end{aligned} \right\}_{\forall km \in \mathcal{E}} 
    \end{align}
\end{subequations}
}

\subsection{DC-OPF Approximation}

Finally, we denote by $\mathcal{Q}^{\text{DC}}$ the DC-OPF model, which is a linear approximation for AC-OPF widely used in practice~\cite{stott+etal09,cain+etal12}. It is based on a number of simplifications that are approximately valid under normal system operations; under this formulation active power losses are zero and reactive power is completely ignored~\cite{bergen+vittal99}.



\subsection{Linear Cuts that Outer-Approximate SOCs}\label{sec:cuts}

In a recent paper~\cite{bichler_2023_prices}, it is empirically shown that the relaxation $\mathcal{Q^{\text{Jabr}}}$ yields less biased price signals than $\mathcal{Q^{\text{DC}}}$ under IP and CH pricing. These findings were derived from small cases in the ARPA-E Grid Optimization Competition 2 (GOC2)~\cite{go2} dataset. It is known that the nonlinear relaxations are simply out of reach for nonlinear solvers on medium-to-large instances~\cite{coffrin_2016_qc, villagra_2024_accurate, bichler_2023_prices}, hence
we explore a {\it linear} relaxation, based on outer-approximating $Q^{\text{Jabr}}$ via linear inequalities. In particular, note that the rotated-cone inequality $x^{2} + y^{2} \leq wz$ is equivalent to $(2x)^{2} + (2y)^{2} \leq (w+z)^{2} - (w-z)^{2}$. Hence,
\begin{equation*}
    x^{2} + y^{2} \leq wz \iff ||(2x,2y,w-z)^{\top}||_{2} \leq w+z. \label{eq:rotatedrewrite}
\end{equation*}
This observation shows that the nonlinear inequalities in the relaxation $\mathcal{Q}^{\text{Jabr}}$ can be written as SOC constraints. Moreover, it can be shown that if some $(x',s') \in \mathbf{R}^{n} \times \mathbf{R}_{+}$ \emph{violates} the SOC inequality $|| x ||_{2} \leq s$, i.e., $||x'||_{2} > s'$, then the hyperplane which achieves maximum Euclidean distance from $(x',s')$ to any hyperplane \emph{separating} the set $\{(x,s) : ||x||_{2} \leq s \}$ from $(x',s')$ is given by $(x')^\top x \leq ||x'||_{2} s$ (see Proposition 5 in~\cite{villagra_2024_accurate}). The latter routine provides a fast procedure for separating over SOC inequalities using purely linear inequalities also known as \emph{linear cuts}. 
We denote by $\mathcal{Q}^{\text{CP}}$ the region resulting from outer-approximating $\mathcal{Q}^{\text{Jabr}}$ by linear cuts (c.f.~\ref{sec:cutplane_algo}). Moreover, given that the SOC inequalities do not depend on input data such as loads, the linear cuts remain valid and can be used if the associated branch remains operational. Hence, we will use the warm-start feature introduced in~\cite{villagra_2024_accurate,villagra_2024_CDC}.

\section{Algorithmic Price Formation} \label{sec:algorithmic-price-formation}
A full description of an electricity price formation algorithm consists of three elements: (1) a model for operations, (2) a parametrization, and (3) a pricing rule \cite{mays_2024_sequential}. We do not discuss parametrization. Instead, we focus on the tradeoff imposed by the AC equations on the first choice, i.e., the model for operations, which is critical for tractability in the third part, i.e., the pricing rule.   

\subsection{Convex Hull and Integer-Programming Pricing}

CH pricing~~\cite{oneil_2005_ipprices} replaces nonconvex cost functions and constraints by their convex hull.
Under piece-wise linear bids, CH prices can be efficiently computed by relaxing the binary constraints \eqref{eq:acopf-welfare-binary-commit} and obtaining the dual variables of the associated convex problem~\cite{hua_2017_convexhull}. IP pricing~\cite{oneil_2005_ipprices} requires solving a mixed integer convex program, fixing the binary variables to their welfare-maximizing value $\{x^{on*}, x^{su*}, x^{sd*} \}$, and obtaining the dual variables from the resulting convex program~\cite{oneil_2005_ipprices}. We refer to \cite{bichler_2023_prices} for a detailed discussion on the properties of these schemes. Critically, both the CH and IP schemes require the underlying model for operations, i.e., the Cartesian product $\mathcal{P}_\mathcal{L} \times \mathcal{P}_\mathcal{G} \times \mathcal{Q}$, to be convex. The DC approximation $\mathcal{Q}^{\text{DC}}$ yields a polyhedral set by giving up representation accuracy. Using $\mathcal{Q}^{\text{Jabr}}$ results in a more precise, yet nonlinear, formulation, which faces scalability issues. Our \emph{Cutting-Plane Pricing Algorithm} provides an intermediate option. 

\subsection{Cutting-Plane Pricing Algorithm}\label{sec:cutplane_algo}

We extend the approach in~\cite{villagra_2024_accurate} to derive accurate electricity prices for large power grids. First, the algorithm constructs a set $\mathcal{Q}^{\text{CP}}$ for~\eqref{model:acopf-welfare-problem} by dynamically generating linear cuts that outer-approximate $\mathcal{Q}^{\text{Jabr}}$. The cuts are implemented under an adequate {\it cut management procedure} which involves (i) quick cut separation; (ii) efficient violated cut selection and (iii) dynamic cut refinement, with online removal of {\it nearly-parallel} and {\it expired} cuts. This yields a rapid and numerically stable algorithm that gives a very tight relaxation with respect to $\mathcal{Q}^{\text{Jabr}}$. For a complete discussion see \cite{villagra_2024_accurate}. 
Since~\eqref{model:acopf-welfare-problem} under $\mathcal{Q}^{\text{CP}}$ is a mixed integer linear program, it can be paired with either CH or IP pricing, see Algorithm \ref{alg:cutting-plane-pricing}. Rel($\mathcal{Q}^{\text{Jabr}}$) refers to the linear inequalities in $\mathcal{Q}^{\text{Jabr}}$, and RelBin($M$) stands for relaxing binary constraints in $M$. 

\setlength{\textfloatsep}{1pt}
\begin{algorithm}
\caption{Cutting-Plane Pricing Algorithm (CPPA)} \label{alg:cutting-plane-pricing}
\begin{algorithmic}[1]
\Procedure {Cutting-Plane Pricing}{}
\State Initialize $r \gets 0$, $z_{0} \gets + \infty$
\State $M \gets \text{Problem \eqref{model:acopf-welfare-problem}}$ with $\mathcal{Q} = \text{Rel}(\mathcal{Q^{\text{Jabr}}})$
\State \textbf{If} CH Pricing \textbf{Then} $M \gets\text{RelBin(M)}$
\While{$t < T$ and $r < T_{ftol}$}
\State $z \gets \min M$ and $\bar{x} \gets \text{argmin} \, M$
\State \textbf{If} $\overline{x} \notin \mathcal{Q}^{\text{Jabr}}$: Compute high quality cuts
\State \hspace{2em} Add cuts if not too parallel to cuts in $M$
\State \hspace{2em} Drop cuts of age $\geq T_{age}$ which are not tight
\State \textbf{If } $z - z_{0} < z_{0} \cdot \epsilon_{ftol}$ \textbf{ Then } $r \gets r+1$
\State \textbf{ElseIf } $z - z_{0} \geq z_{0} \cdot \epsilon_{ftol}$ \textbf{ Then } $r \gets 0$
\State $z_{0} \gets z$
\EndWhile
\State Solve $M$
\State \textbf{If} IP Pricing \textbf{Then} Fix binary vars in $M$, re-solve
\State Extract dual variables as prices $\lambda$
\EndProcedure
\end{algorithmic} 
\end{algorithm} 

Relevant input parameters for our procedure are: a time limit $T>0$; a number of admissible iterations without sufficient objective improvement $T_{ftol} \in \mathbb{N}$; and a threshold for relative objective improvement $\epsilon_{ftol} > 0$. Our measure of {\it cut-quality} is the amount by which a solution $\bar{x}$ to $M$ violates some nonlinear inequality in $\mathcal{Q}^{\text{Jabr}}$, among which we add a top percentage to $M$ (c.f.~\ref{sec:cuts}). We filter cuts by age ($T_{age}$) and a measure of {\it parallelism} - see \cite{villagra_2024_accurate} for more details.

\section{Numerical Results} \label{sec:numerical-results}

\setlength{\tabcolsep}{4pt} 
\begin{table*}[t]
\caption{Price Statistics and Run Time for Jabr SOC, DC and CP with IP}
\centering
\begin{tabular}{ @{} l r r r r r r r r r r r r @{} }
\toprule
& \multicolumn{3}{c}{Jabr SOC Pricing} & \multicolumn{4}{c}{DC Pricing} & \multicolumn{5}{c}{CP Pricing} \\
\cmidrule(l{0.5em}r{0.40em}){2-4} \cmidrule(l{0.5em}r{0.40em}){5-8} \cmidrule(l{0.5em}r{0.40em}){9-13}
\multicolumn{1}{l}{Scenario} & Avg & Std & Time & Avg & Std & $\delta$ & Time & Avg & Std & $\delta$ & Time & \#Cuts \\
\midrule
617-005 & 95.67 & 1.17 & 13.44 & 94.81 & 1.4{\it e}-11 & 1.15 & 0.97 & 95.55 & 1.16 & 0.11 & 9.61 & 1,945 \\
617-017 & 118.64 & 1.24 & 43.46 & 118.74 & 5.5{\it e}-12 & 0.84 & 0.94 &  118.58 & 1.23 & 0.07 & 12.07 & 2,052 \\
617-024 & 61.32 & 1.10 & 5.18 & 59.94 & 4.8{\it e}-11 & 1.50  & 1.07 &  61.26 & 1.09 & 0.07 & 7.41 & 1,917\\
617-062 & 82.26 & 1.25 & 14.22 & 79.06 & 3.7{\it e}-11 & 3.24 & 1.02 &  82.23 & 1.24 & 0.04 & 7.35 & 1,558 \\
617-073 & 106.61 & 1.34 & 10.47 & 106.48 & 1.1{\it e}-12 & 0.91 & 1.12 &  106.53 & 1.34 & 0.08 & 7.43 & 1,963 \\
\midrule
617-Avg & 92.90 & 1.22 & 17.35 & 91.81 & 2.1{\it e}-11 & 1.53 & 1.03 &  92.83 & 1.21 & 0.07 & 8.77 & 1,887  \\
\bottomrule
\end{tabular}
\label{table:617-prices-&-runtimes}
\vspace{-5pt}
\end{table*}

\setlength{\tabcolsep}{3pt} 
\begin{table*}[t]
\caption{Economic Efficiency Metrics for Jabr SOC, DC and CP Relaxations with IP after Redispatch}
\centering
\begin{tabular}{ @{} l r r r r r r r r r r r r r r r @{} }
\toprule
& \multicolumn{5}{c}{Jabr SOC Pricing} & \multicolumn{5}{c}{DC Pricing} & \multicolumn{5}{c}{CP Pricing} \\
\cmidrule(l{0.5em}r{0.40em}){2-6} \cmidrule(l{0.5em}r{0.40em}){7-11} \cmidrule(l{0.5em}r{0.40em}){12-16}
\multicolumn{1}{l}{Scenario} & Welfare & MWPs & GLOCs & LLOCs & RDCs & Welfare & MWPs & GLOCs & LLOCs & RDCs & Welfare & MWPs & GLOCs & LLOCs & RDCs\\
\midrule
617-005 & 646,240 & 37,429 & 230,322 & 196,674 & 331 & 645,158 & 36,974 & 231,608 & 198,404 & 6,866 & 646,226 & 37,432 & 230,357 & 196,709 & 635 \\
617-017 & 573,962 & 47,221 & 299,280 & 254,546 & 335 & 573,447 & 47,209 & 299,824 & 255,091 & 5,915 & 573,967 & 47,222 & 299,287 & 254,553 & 1,235 \\
617-024 & 741,538 & 9,739 & 141,700 & 138,485 & 288 & 739,697 & 9,749 & 143,844 & 140,628 & 6,876 & 741,533 & 9,742 & 141,719 & 138,504 & 468 \\
617-062 & 701,515 & 19,339 & 177,561 & 163,122 & 318 & 700,209 & 19,397 & 179,414 & 164,975 & 7,035 & 701,508 & 19,340 & 177,572 & 163,133 & 560 \\
617-073 & 610,392 & 48,218 & 265,102 & 220,130 & 366 & 609,556 & 48,202 & 266,001 & 221,029 & 6,844 & 610,388 & 48,220 & 265,119 & 220,146 & 643 \\
\midrule
617-Avg & 654,729 & 32,389 & 222,793 & 194,591 & 327 & 653,631 & 32,306 & 224,138 & 196,026 & 6,707 & 654,724 & 32,391 & 222,811 & 194,609 & 712 \\
\bottomrule
\end{tabular}
\label{table:617-welfare-redispatch-costs}
\vspace{-5pt}
\end{table*}

\subsection{Experimental Setup}

We report numerical experiments on scenarios for a 617-bus instance and the three largest instances from the GOC2 dataset~\cite{go2}. As in~\cite{bichler_2023_prices}, we abstract from contingencies, switched shunts, and line and transformer switching. All of our experiments were run on an Intel(R) Xeon(R) Linux64 machine CPU E5-2687W v3 3.10 GHZ with 20 physical cores, 40 logical processors and 256 GB RAM. We used three commercial solvers: Gurobi version 10.0.1~\cite{gurobi}, Artelys Knitro version 13.2.0~\cite{knitro}, and Mosek 10.0.43~\cite{mosek}. Gurobi was used to solve the LPs for CPPA, Mosek was used to solve the Jabr SOCs, and Knitro for finding AC-feasible solutions. Our code and solution files can be downloaded from www.github.com/matias-vm.

For the first case study (617), we pair the Jabr SOC, DC and CP formulations with the IP pricing rule to obtain an allocation  $\allocation^{\relaxation}=(v,p,x)$, and prices $\lambda^\relaxation$ for $\relaxation$ in $[\text{Jabr}, \text{DC}, \text{CP}]$. First, we directly compare prices and computation times for all the above models (Table \ref{table:617-prices-&-runtimes}). We then follow \cite{bichler_2023_prices} and simulate a {\it redispatch} procedure by taking the optimal allocation $\allocation^{\relaxation}$ and feeding it as a warm-start to Knitro, to find an adjusted AC-feasible allocation, denoted $\phi(\allocation^\relaxation)$. 
Using this allocation, as in~\cite{bichler_2023_prices}, we compute metrics of economic efficiency (Table \ref{table:617-welfare-redispatch-costs}) (i) make-whole payments (MWPs), (ii) global and local lost opportunity costs (GLOCs/LLOCs), and (iii) re-dispatch costs (RDCs) which we compare against welfare. Formally, for an allocation $\allocation$ and prices $\lambda$, consider the \textit{direct utility} as $\utility_g(\allocation,\lambda) = p_g\lambda_{k(g)}- (-\pi_g(\allocation))$ or $\utility_l(\allocation,\lambda) = \pi_l(\allocation) - p_l\lambda_{k(l)}$ where $k(\tau)$ for $\tau\in\{l,g\}$ indicates the bus associated to the load or generator. The metrics are then defined as

\begin{subequations}
{\small
\begin{align}
& MWP  := \sum_{\tau \in \mathcal{L}\cup \mathcal{G}} [- \utility_\tau(\phi(\allocation),\lambda) ]^+
 \\
& GLOC  := \sum_{\tau \in \mathcal{L}\cup \mathcal{G}} [ \sup_{\allocation'} \{\utility_\tau(\allocation',\lambda)\} - \utility_\tau(\phi(\allocation),\lambda)] \\
& LLOC  := \sum_{\tau \in \mathcal{L}\cup \mathcal{G}} [ \sup_{p_\tau} \{\utility_\tau(p_\tau,p_{-\tau},x,\lambda)\} - \utility_\tau(\phi(\allocation),\lambda)]\\
& RDC := \sum_{\tau \in \mathcal{L}\cup \mathcal{G}} [ \utility_\tau(\allocation,\lambda) -\utility_\tau(\phi(\allocation),\lambda) ]^+
\end{align}
}
\end{subequations}

For the larger instances, we pair the Jabr SOC, DC and CP formulations with the CH pricing rule. We directly compare prices and run times for all the above models (c.f. Table~\ref{table:price-statistics-in-large-grids}). In Table~\ref{table:price-statistics-in-large-grids-contingencies} we report on scenarios under contingencies.



\setlength{\tabcolsep}{3pt} 
\begin{table*}[t]
\caption{Price Statistics and Run Time for Jabr SOC, DC and CP with CH Policy}
\centering
\begin{tabular}{ @{} l r r r r r r r r r r r r r r r r @{} }
\toprule
& \multicolumn{4}{c}{Jabr SOC Pricing} & \multicolumn{5}{c}{DC Pricing} & \multicolumn{7}{c}{CP Pricing} \\
\cmidrule(l{0.5em}r{0.40em}){2-5} \cmidrule(l{0.5em}r{0.40em}){6-10} \cmidrule(l{0.5em}r{0.40em}){11-17}
\multicolumn{1}{l}{Scenario} & Avg & Std & Time & Obj & Avg & Std & $\delta$ & Time & Gap\% & Avg & Std & $\delta$ & TimeLP  & TimeCut & \#Cuts & Gap\% \\
\midrule
17,700-019 & 1.01 & 1.07 & 66.06 & 7,789,436 & 15.12 & 0.61 & 13.38 & 19.73 & 9.40 & 14.82 & 1.05 & 13.09 & 19.33 & 285.40 & 69,364 & 9.659 \\
17,700-020 & 1.22 & 1.62 & 73.55 & 7,757,103 & 14.34 & 1.12 & 12.45 & 21.96 & 9.260 &  14.72 & 1.30 & 12.80 & 23.88 & 285.48 & 70,744 & 9.535  \\
17,700-021 & 1.26 & 1.58 & 70.68 & 7,705,752 & 13.59 & 1.21 & 11.69 & 26.43 & 8.874 & 14.08 & 1.33 & 12.15 & 34.67 & 290.93 & 71,346 & 9.163  \\
17,700-089 & 0.80 & 0.71 & 64.22 & 7,792,189 & 14.71 & 0.73 & 13.18 & 24.43 & 9.290 & 14.71 & 1.06 & 13.20 & 22.83 & 304.61 & 70,043 & 9.518 \\
17,700-094 & 1.09 & 1.25 & 66.70 & 7,778,704 & 13.99 & 1.01 & 12.24 & 20.59 & 8.390 &  14.28 & 1.21 & 12.52 & 21.62 & 287.68 & 70,223 & 8.659 \\
\midrule
17,700-Avg & 1.08 & 1.25 & 68.24 & 7,764,637 & 14.35 & 0.94 & 12.59 & 22.63 & 9.043 & 14.52 & 1.19 & 12.75 & 24.47 & 290.82 & 70,344 & 9.307  \\
\midrule
19,402-006 & 1.69 & 0.036 & 150.22 & 233,510 & 1.66 & 0.078 & 0.049 & 28.24 & 0.92 & 1.68 & 0.04 & 0.002 & 48.30 & 331.12 & 147,868 & 0.053 \\
19,402-010 & 1.33 & 0.030 & 153.81 & 280,531 & 1.32 & 0.11 & 0.055 & 28.31 & 0.45 & 1.33 & 0.03 & 0.002 & 48.13 & 345.71 & 147,857 & 0.047  \\
19,402-069 & 1.64 & 0.038 & 155.89 & 478,515 & 1.62 & 0.084 & 0.052 & 26.09 & 0.36 & 1.64 & 0.04 & 0.003 & 44.16 & 307.01 & 143,397 & 0.060  \\
19,402-077 & 1.35 & 0.022 & 137.88 & 202,145 & 1.34 & 0.043 & 0.029 & 27.48 & 0.67 & 1.35 & 0.02 & 0.003 & 45.89 & 323.92 & 143,657 & 0.095 \\
19,402-095 & 1.48 & 0.023 & 144.32 & 225,801 & 1.47 & 0.027 & 0.024 & 27.00 & 0.57 & 1.47 & 0.02 & 0.003 & 51.19 & 331.81 & 141,628 & 0.031  \\
\midrule
19,402-Avg & 1.50 & 0.030 & 148.42 & 284,100 & 1.48 & 0.069 & 0.042 & 27.43 & 0.60 & 1.49 & 0.03 & 0.002 & 47.54 & 327.91 & 144,881 & 0.057   \\
\midrule
31,777-011 & 3.95 & 0.92 & 146.40 & 64,740,242 & 4.53 & 11.14 & 1.15 & 37.13 & 0.14 & 3.91 & 0.34 & 0.039 & 56.51 & 299.09 & 137,764 & 0.002 \\
31,777-012 & 3.10 & 0.30 & 127.06 & 52,756,414 & 3.55 & 7.92 & 1.08 & 36.35 & 0.12 & 3.09 & 0.30 & 0.012 & 40.39 & 262.79 & 137,764 & 0.002  \\
31,777-015 & 3.36 & 1.82 & 151.64 & 62,000,434 & 3.92 & 12.27 & 1.39 & 36.45 & 0.15 & 3.34 & 1.49 & 0.026 & 40.03 & 262.42 & 136,885 & 0.002 \\
31,777-019 & 3.12 & 0.87 & 139.62 & 55,916,108 & 3.61 & 8.80 & 1.15 & 37.06 & 0.11 & 3.10 & 0.48 & 0.039 & 40.19 & 275.03 & 136,536 & 0.002  \\
31,777-051 & 3.42 & 0.29 & 132.43 & 56,825,370 & 4.17 & 11.24 & 1.40 & 36.01 & 0.15 & 3.41 & 0.29 & 0.013 & 33.66 & 298.17 & 120,857 & 0.002  \\
\midrule
31,777-Avg & 3.39 & 0.84 & 139.43 & 58,447,713 & 3.96 & 10.27 & 1.23 & 36.60 & 0.14 & 3.37 & 0.58 & 0.026 & 42.16 & 279.50 & 133,991 & 0.002  \\
\bottomrule
\end{tabular}
\label{table:price-statistics-in-large-grids}
\vspace{-10pt}
\end{table*}

\setlength{\tabcolsep}{3pt} 
\begin{table*}[t]
\caption{Price Statistics and Run Time for Jabr SOC, DC and CP with CH Policy under Contingencies}
\centering
\begin{tabular}{ @{} l r r r r r r r r r r r r r r r r @{} }
\toprule
& \multicolumn{4}{c}{Jabr SOC Pricing} & \multicolumn{5}{c}{DC Pricing} & \multicolumn{7}{c}{CP Pricing} \\
\cmidrule(l{0.5em}r{0.40em}){2-5} \cmidrule(l{0.5em}r{0.40em}){6-10} \cmidrule(l{0.5em}r{0.40em}){11-17}
\multicolumn{1}{l}{Scenario} & Avg & Std & Time & Obj & Avg & Std & $\delta$ & Time & Gap\% & Avg & Std & $\delta$ & TimeLP  & TimeCut & \#Cuts & Gap\% \\
\midrule
31,777-011 (C) & * & - & 162.63 & - & 5.09 & 13.14 & - & 36.47 & - & Infeas & - & - & 78.97 & - & 129,689 & - \\
31,777-012 (C) & * & - & 138.43 & - & 3.55 & 7.92 & - & 35.93 & - & Infeas & - & - & 65.69 & - & 107,467 & - \\
31,777-015 (C) & 3.45 & 3.66 & 149.27 & 61,985,332 & 4.33 & 12.20 & 1.78  & 36.23 & 0.16 & 3.42 & 3.60 & 0.029 & 72.26 & - & 108,718 & 0.150 \\
31,777-019 (C) & * & - & 103.61 & - & 3.94 & 10.26 & - & 36.36 & - & Infeas & - & - & 80.05 & - & 107,606 & - \\
31,777-051 (C) & 3.50 & 2.79 & 147.85 & 56,810,808 & 4.26 & 11.95 & 1.61 & 36.44 & 0.176 & 3.47 & 2.41 & 0.039 & 37.56  & - & 109,477 & 0.003 \\
\bottomrule
\end{tabular}
\label{table:price-statistics-in-large-grids-contingencies}
\vspace{-10pt}
\end{table*}

\subsection{Computational Results}\label{subsection:experiments}



\begin{figure}
    \centering
    \caption{Average price distance $\delta$ for CP and DC to Jabr SOC}
    \includegraphics[width=0.98\linewidth]{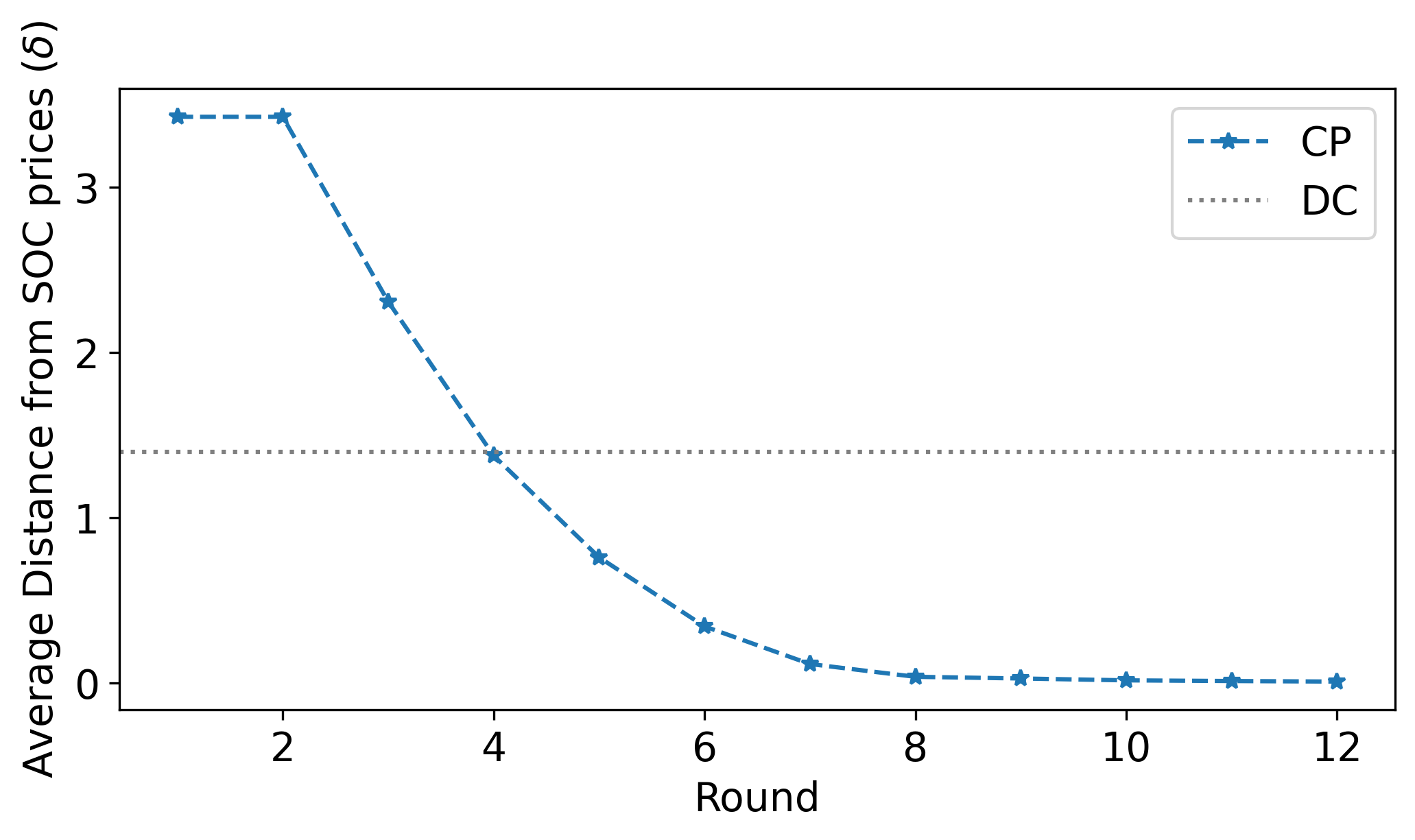}
    \label{fig:price_diff_solo}
\end{figure}


\paragraph{Case 617} 
Table \ref{table:617-prices-&-runtimes} displays the average (``Avg") and standard deviation (``Std") of nodal prices (\$/MWh) for each scenario, as well as the average across scenarios, for the Jabr SOC, the DC approximation, and the CP relaxation.\footnote{We conducted a thorough proof check of our experimental setting and code while deriving our results. Our code is available online. Still, we acknowledge small discrepancies with results obtained for the 617-bus instances with \cite{bichler_2023_prices}.} DC produces a slightly lower average price across nodes and scenarios 
, while attaining zero variance across nodes failing to capture network congestion. Conversely, the SOC and CP relaxations yield prices with non-zero --- and very similar --- values for standard deviation. Moreover, CP prices are extremely close to those arising from Jabr SOCP as measured by average distance (``$\delta$"), formally defined as $\delta({\text{rel}}) = \|\lambda_{\text{rel}}-\lambda_{\text{SOC}}\|_1/|\mathcal{B}|$ for each relaxation ``rel" and network $\mathcal{N}=(\mathcal{B},\mathcal{E})$. 
The run times (``Time") for each model in Table~\ref{table:617-prices-&-runtimes} show the run time per scenario in seconds, as well as the average across scenarios. The number of cuts (\#Cuts) generated for each scenario is shown in the right-most column of the table.




Table~\ref{table:617-welfare-redispatch-costs} displays economic efficiency metrics (in \$) with respect to the allocation and prices produced by the Jabr SOC, DC, and CP formulations {\it after the redispatch process}. 
After AC-feasibility adjustments, welfare is slightly higher for the Jabr SOC and CP relaxations than for the DC approximation. More importantly, {\it redispatch costs} are one order of magnitude smaller: in the case of DC, the RDCs are, on average, a $1.02 \%$ of total welfare, while Jabr SOC and CPPA show relative values of $0.05 \%$ and $0.11 \%$ instead. 



\paragraph{Large cases ($\geq 15,000$ nodes)}

Table~\ref{table:price-statistics-in-large-grids} displays the average and standard deviation of nodal prices for all the considered scenarios. We ran the CPPA for $300$ seconds and report the total time spent generating cuts (``TCut"), as well as the time it takes to solve the last LP (``TimeLP"). We report the objective (in \$) of~\eqref{model:acopf-welfare-problem} for Jabr SOC (``Obj") and the \% gap relative to DC and CP objectives (``Gap\%"). Overall, CPPA achieves very quickly prices that are close (small $\delta$) to those attained by Jabr SOC; along the lines of what is reported in~\cite{bichler_2023_prices} for small networks. The run times of the last LPs are comparable to the those of the DC model and faster than Jabr SOC, which serves as an estimate of how long CPPA would take if warm-started. For the 17,700-bus scenarios, we find that neither DC nor CP resemble Jabr SOC prices. Further investigation of these scenarios let us confirm that these instances are actually AC-infeasible, via a computation of a minimally infeasible system for the i2 relaxation for AC-OPF~\cite{farivar_2011_i2,coffrin_2016_qc,villagra_2024_accurate}, which is stronger than the Jabr SOC. Hence, in these scenarios, Jabr SOC fails to capture the physics of the AC power flows.






\paragraph{Warm-Starts for 31,777-bus}

We tested the warm-starting capabilities for CPPA by considerably stressing the five scenarios of the 31,777-bus network (see~Table~\ref{table:price-statistics-in-large-grids-contingencies}). To this end, first we ran CPPA for 400 seconds and saved the cuts computed in the last round. We observe in Figure~\ref{fig:price_diff_solo} the improvement, in terms of average distance $\delta$, of every round of the algorithm in scenario 51. We created stressed instances by loading branch contingencies, at most 1\% of the total branches of the network, provided in the GOC2 dataset for each scenario. Next, we warm-started CPPA with the previously computed cuts (``\#Cuts"), and run CPPA for only one round. Mosek failed to converge for scenarios 11, 12 and 19 under contingencies (``$*$"), whereas CPPA quickly identified these scenarios as infeasible (``Infeas"), and DC failed to do so. In fact, since these scenarios are infeasible, the price signals reported by the DC model are entirely inaccurate, with average prices being significantly underestimated. For scenarios 15 and 51, CPPA attains, twice as fast, prices that are very close to Jabr SOC, with $\delta(\text{CP})$ values being at least an order of magnitude smaller than $\delta(\text{DC})$ values. In Table~\ref{table:617-reactive-power-prices} we display reactive power prices for the SOC and CP relaxations. Note that the average reactive power prices of the CP relaxation are very close to the Jabr SOC, however, there is a larger nodal price variation. The DC approximation does not produce reactive power prices.

\vspace{-0.5em}

\setlength{\tabcolsep}{3pt} 
\begin{table}[h!]
\caption{Reactive Power Prices for Jabr SOC, DC and CP with CH}
\centering
\begin{tabular}{ @{} l r  r r r r r r r r @{} }
\toprule
& \multicolumn{2}{c}{Jabr SOC Pricing} & \multicolumn{2}{c}{DC Pricing} & \multicolumn{3}{c}{CP Pricing} \\
\cmidrule(l{0.5em}r{0.40em}){2-3} \cmidrule(l{0.5em}r{0.40em}){4-5} \cmidrule(l{0.5em}r{0.40em}){6-8}
\multicolumn{1}{l}{Scenario} & Avg & Std & Avg & Std & Avg & Std & $\delta$ \\
\midrule
31,777-011 & -0.0069 & 1.76{\it e}-5 & - & - & -0.0079 & 0.67 & 0.14 \\
31,777-012 & -0.0159 & 2.02{\it e}-5 & - & - & -0.0163 & 0.06 & 0.02 \\
31,777-015 & 0.0168 & 1.6{\it e}-4 & - & - & 0.0269 & 1.32 & 0.60 \\
31,777-019 & -0.0120 & 8.2{\it e}-6 & - & - & -0.0126 & 0.15 & 0.05 \\
31,777-051 & -0.0113 & 8.2{\it e}-7 & - & - & -0.0124 & 0.08 & 0.10 \\
\midrule
31,777-Avg & -0.0059 & 4.1{\it e}-5 & - & - & -0.0045 & 0.46 & 0.24 \\
\bottomrule
\end{tabular}
\label{table:617-reactive-power-prices}
\vspace{-5pt}
\end{table}

\section{Conclusions} \label{sec:conclusions-future-work}

We develop a pricing algorithm based on a linear relaxation of the AC-OPF problem that tightly approximates the Jabr SOC. We conducted numerical experiments to measure accuracy of price signals, redispatch costs, and computational scalability. We show that the our algorithm, using IP pricing, delivers highly accurate prices and low redispatch costs compared to DC on a 617-bus system. For the three largest grids in GOC2, our algorithm CPPA with CH pricing achieves prices close to the Jabr SOC with run times comparable to the DC approximation when warm-started.




 \section{Acknowledgements}

\thanks{We would like to thank Daniel Bienstock, Richard O'Neill, Johannes Knorr and participants at the INFORMS 2024 Annual Meeting for insightful comments and discussions.}



\end{document}